\documentclass[12pt]{article}
\usepackage{amsfonts, amsmath, amssymb, amsgen, amsthm, amscd,
latexsym}
\usepackage[all]{xy}

\def\Bc{{\cal B}}
\def\Cc{{\cal C}}
\def\Dc{{\cal D}}

\def\Hc{{\cal H}}

\def\Cc{{\cal C}}

\def\a{\alpha}

\def\D{\Delta}

\def\0b{\bf 0}

 \def\bbc{\mathbb{C}}

\def\ot{\otimes}

\def\ra{\rightarrow}

\def\rt{\rtimes}

\def\0D{\Delta^{(0)}}
\def\1D{\Delta^{(1)}}

\newtheorem{theorem}{Theorem}[section]

\newtheorem{proposition}[theorem]{Proposition}
\newtheorem{lemma}[theorem]{Lemma}
\newtheorem{corollary}[theorem]{Corollary}

\def\build#1_#2^#3{\mathrel{
\mathop{\kern 0pt#1}\limits_{#2}^{#3}}}
\newcommand{\ps}[1]{~\hspace{-4pt}_{^{(#1)}}}

\numberwithin{equation}{section}

\begin{document}

\title{\bf Constant  and Equivariant Cyclic Cohomology}

\author{  Bahram Rangipour   \\ \\
    Department of Mathematics \\
    The Ohio State University \\
    Columbus, OH 43210, USA
}
\date{ \ }

\maketitle

\begin{abstract}
In this note we prove that  the constant and equivariant cyclic
cohomology of algebras coincide. This shows that constant cyclic
cohomology is rich and computable.
\end{abstract}

\section{Introduction}
Cyclic cohomology, invented  by A. Connes \cite{ac}  as a
cornerstone of noncommutative geometry,  is the noncommutative
analogue of de Rahm homology. It plays a crucial role  in
noncommutative geometry by detecting topological invariants such as
$K$-theory of algebras via a pairing  which is a  nontrivial
generalization of  the index of operators in the classical geometry.
There is a variety of specialized
 cyclic cohomology theories,  which are defined in various contexts
 for different purposes. In this
 paper we briefly discuss two of these  theories and prove that they are
 actually the same.

 Symmetry in classical geometry is defined via  the actions of groups on
  spaces. In noncommutative geometry  spaces  are  replaced
 with  not necessarily commutative algebras, viewed as algebras of coordinates
 of "noncommutative spaces". Accordingly,  the role of symmetry is played by
 actions of  Hopf algebras on algebras. A natural entry into
 this page of the `commutative-noncommutative' dictionary  is the counterpart
of equivariant de Rahm homology, which evidently is the
 equivariant cyclic cohomology. Equivariant cyclic cohomology has been studied by
   various authors, in particular we refer the reader to \cite{b,bg,ft, gj, n1,n2},
    for the case of groups, and \cite{ak,nt,hkrs2}, for Hopf algebras.

In homological algebra, it is often important   to replace a complex
with a quasisomorphic subcomplex. In cyclic cohomology it is known
that normalized co-chains form a subcomplex quasisomorphic to cyclic
complex.  The former  are those co-chains that vanish on the image
of degeneracies. Restricting on the normalized subcomplex is usually
useful and sometime even necessary, for example   in the proof of
the cyclic Eilenberg-Zilber theorem. Connes in \cite{ac04} in the
study of spectral triple over $A(SU_q(2))$ generalized the
normalized cyclic complex to constant cyclic subcomplex. Let us
briefly recall it. Let $\Cc$ be a unital subalgebra of a unital
algebra $\Dc$. A $\Cc$-constant co-chain is a  co-chain $\varphi\in
{\rm Hom}(D^{\ot (n+1)},\bbc)$, such that $\varphi(a^0, \dots,
a^n)$, and $b\varphi(a^0, \dots, a^{n+1})$ vanish whenever $\a^j\in
\Cc$, for at least one  $j\ge 1$. Then the $\Cc$-constant co-chains
form a sub-complex of the mixed-complex computing the cyclic
cohomology of $\Dc$ and hence one defines $\Cc$-constant cyclic
cohomology of $\Dc$ as the cohomology of this subcomplex.

Now let us consider a Hopf algebra $\Hc$, with invertible antipode,
acting on an algebra $A$ and makes $A$ an $\Hc$-module algebra. One
knows that the algebra $\Hc$ is a subalgebra of the crossed product
algebra $A\rt \Hc$. So we can talk about $\Hc$-constat cohomology of
$A\rt \Hc$. On the other hand equivariant cyclic cohomology of $A$
under the action of $\Hc$ is well-defined. In this paper we prove
that  the two theories coincide.

Throughout this note we assume that algebras are unital and
associative,  that every algebraic structure including Hopf algebras
are defined over
 $\bbc$, and  that Hopf algebras have invertible
antipode. The antipode, coproduct and counit of Hopf algebras are
denoted by $S$, $\D$, and $\epsilon$ respectively. The action of
$\D$ on an element of $\Hc$, say $h$, is denoted by $h\ps{1}\ot
h\ps{2}$.

\section{Constant Cyclic Cohomology}
Let $A$ be a unital algebra. We recall, from \cite{ac},  the cyclic
cohomology of $A$ is the cohomology of the space of all  cyclic
co-chains, which are
 those Hochschild co-chains that  are cyclic, i.e, $f\in {\rm Hom}(A^{\ot
(n+1)}, \bbc)$ satisfy $f(a^0,\dots, a^n)=(-1)^n f(a^n,a^0 \dots,
a^{n-1})$. The space of all cyclic co-chains forms a subcomplex of
Hochschild complex of the algebra $A$ with coefficients in $A^\ast$,
the dual space  of $A$.

Cyclic cohomology of $A$ can be also defined via the Connes
$(b,B)$-bicomplex. Let us  recall  this bicomplex. let $\Bc^{n,m}=
\Bc^{n,m}(A)= A^{\ot(n-m)}$ and $0$ where $m\; >\; n$ or $m\;<\;0$.
With the Connes boundary map $B: \Bc^{n,m}\ra \Bc^{n,m+1}$, and the
Hochschild's one $b: \Bc^{n,m}\ra \Bc^{n+1,m}$, $\Bc$ forms a
bicomplex.  The cyclic complex  of $A$ is quasisomorphic  to  the
total complex of $\Bc(A)$. The vertical and horizontal boundary maps
are defined as follows.

First $B:=AB_0$, where \begin{equation} (B_0 \varphi)(a^0,\dots,
a^n):=\varphi(1, a^0,a^1,\dots, a^n)-(-1)^{n+1}\varphi(a^0,a^1\dots,
a^n,1)
\end{equation}
\begin{equation}
A\varphi(a^0,a^1,
\dots,a^n):=\sum_{0}^n(-1)^{nj}\varphi(a^j,a^{j+1}\dots,a^{j-1}),
\end{equation}
and the Hochschild boundary map,
\begin{align}
&b\varphi(a^0, \dots, a^{n+1}):=\\
&\sum_0^n (-1)^j\varphi(a^0, \dots, a^ja^{j+1}, \dots,
a^{n+1})+(-1)^{n+1}\varphi(a^{n+1}a^0, a^1,\dots, a^{n}).
\end{align}

 Let us recall the definition
of constant cyclic cohomology from \cite{ac04}. Let $\Cc\subset \Dc$
be a unital subalgebra of a unital algebra $\Dc$. We say a co-chain
$\varphi\in {\rm Hom}(\Dc^{\ot (n+1)}, \bbc)$ is  $\Cc$-constant if
and only if both $\varphi(a^0,\dots a^n)$, and $b\varphi(a^0,\dots
a^{n+1})$ vanish whenever  at least one of the $a^j$, $j\ge 1$, is
in $\Cc$.

Since the subalgebra $\Cc$ is unital one sees that
$B_0\varphi(a^0,\dots, a^n)=\varphi(1, a^0, \dots, a^n)$,  for a
$\Cc$-constant co-chain $\varphi$ . This proves that the boundary
$B$ sends $\Cc$-constant co-chains to themselves. So $\Cc$-constant
co-chains form a subcomplex of the $(b,B)$ bicomplex and hence we
can define the $C$-constant cyclic cohomology to be the total
cohomology of this subcomplex. We denote this cohomology by
$HC^*(\Dc;\Cc)$.

If one takes the ground field as the subalgebra $\Cc$ then the
$\Cc$-constat co-chains are just normalized co-chains and hence
$HC^*(\Dc; \bbc)=HC^*(\Dc)$.

Now let $\Hc$ be a Hopf algebra acting on an algebra $A$ and makes
$A$ a $\Hc$-module algebra. The latter means that the multiplication
and unit map of $A$ are $\Hc$-linear, i.e,
$h(ab)=h\ps{1}(a)h\ps{2}(b)$, and $h(1)=\epsilon(h)$. Having the
above situation, we make a crossed product  algebra $A\rt \Hc$,
which has $A\ot \Hc$ as its underlying vector space, $1\ot 1$ as its
unit,  and  $(a\ot h)(b\ot g)=ah\ps{1}(a)\ot h\ps{2}g$\; as its
multiplication. One can see that $\Hc$ is a unital subalgebra of
$A\rt \Hc$, and hence $\Hc$-constant cyclic cohomology of $A\rt \Hc$
is well-defined.

Our main goal in this paper is to compute this cohomology in terms
of equivariant cyclic cohomology defined in \cite{ak,nt} which we
recall it in the next section.
%%%%%%%%%%%%%%%%%%%%
%%%%%%%%%%%%%%%%%%%%%%%%%%%%%%%%%%%
%%%%%%%%%%%%%%%%%%%%%%%%%%%%%%%%%%%%%%%%%%%%%%%%%%%%%%
\section{Equivariant Cyclic Cohomology}
Cyclic cohomology in its ultimate extent  is defined for any
co-cyclic module \cite{ac,cncg}. A co-cyclic module is a
co-simplicial module $M=(\{M^n\}_{n\ge 0}, d_n^i, s_n^j)$, $0\le
i\le n+1$, $0\le j\le n-1$,  together with a cyclic operator
$t_n:M^n\ra M^n$, for each $n=0,1,\dots$ satisfy the following
identities.
\begin{align}
&t_{n+1}d_n^i= d_n^{i-1}t_n,\quad, t_{n+1}d_n^0=d_n, \\
&t_{n-1}s_n^j=s_n^{j-1}t_n, \quad t_{n-1}s_n^0=s_n^{n-1}t_n^2,\\
& t_n^{n+1}=id
\end{align}
The cyclic cohomology of a co-cyclic module is defined to be the
cohomology of total complex of $\Bc(M)$ which is a bicomplex that in
degree $(n,m)$ is $M^{n-m}$  and $0$ above the main diagonal and
below the horizontal axis. Its vertical and horizontal  boundaries
$B$ and $b$ are defined  as follows.

\begin{equation}
b:=\sum_0^{n+1}(-1)^id_n^i,
\end{equation}
\begin{align}
&B:=AB_0, \\
&B_0:= (1+(-1)^{n}t_n)s_n,\\
&A:=\sum_0^n (-1)^{nj}t_n^j
\end{align}
%%%%%%%%%%%%%%%%%%%%%%%

Now, let us recall the definition of equivariant cyclic cohomology.
Equivariant cyclic cohomology of an algebra under the action of a
discrete group is studied  in \cite{b,bg,ft, gj, n1,n2}. Its
generalization for the action of a Hopf algebra on an algebra is
dealt with  in \cite{ak,nt,hkrs2}. In the following we recall this
cohomology theory.

Let
 $C^n_\Hc(A)= {\rm Hom}_\Hc(\Hc\ot A^{\ot (n+1)}, \bbc)$,
 which is all equivariant co-chains, i.e., all
  $\varphi:\Hc\ot A^{\ot (n+1)}\ra \bbc $, such that  for
  all $g,h\in \Hc$, and $a_i\in A$,
  \begin{equation}\label{equiv}
  \varphi(g\ps{1}\cdot h, g\ps{2}(a^0), \dots ,
 g\ps{n+1}(a^n))=\epsilon(g)\varphi(h,a^0, \dots, a^n),
 \end{equation}
 where $g\cdot h= g\ps{1}hS^{-1}(g\ps{2}) $ is the usual adjoint
 action of $\Hc$ on itself.
 \begin{lemma}\label{condition}
The condition \eqref{equiv} is equivalent to
\begin{equation}
f(S(g\ps{2})h g\ps{1}, a^0, \dots, a^n)= f(h, g\ps{1}(a^0), \dots,
g\ps{n+1}(a^n)).
\end{equation}
\end{lemma}
\begin{proof}
 Let $f$ satisfy  the condition \eqref{equiv}.
\begin{align*}
&f(h , g\ps{1}(a^0), \dots, g\ps{n+1}(a^n))=\\
&\epsilon(g\ps{1})f(h , g\ps{2}(a^0), \dots, g\ps{n+2}(a^n))=\\
&\epsilon(S(g\ps{1}))f(h , g\ps{2}(a^0), \dots, g\ps{n+2}(a^n))=\\
&f(S(g\ps{1})\ps{1}\cdot h , S(g\ps{1})\ps{2}g\ps{2}(a^0), \dots, S(g\ps{1})\ps{n+2}g\ps{n+2}(a^n))=\\
&f(S(g\ps{2})h g\ps{1}, a^0, \dots, a^n).
\end{align*}

The converse is similarly checked.
\end{proof}
 One can see that, see one of  \cite{ak,nt, hkrs2},  $C^\ast_\Hc(A)$ with the following operators
 is
 a co-cyclic module.
\begin{align}
&d_n^i\varphi(h,a^0,\dots,a^{n+1})=\varphi(h,a^0,\dots, a^ia^{i+1},
\dots, a^{n+1}), \quad, 0\le i\le n\\
&d_n^{n+1}\varphi(h,a^0, \dots, a^{n+1})= \varphi(h\ps{2},
S^{-1}(h\ps{1})(a^{n+1})a^0,a^1, \dots, a^{n}),\\
&s_n^j\varphi(h, a^0, \dots, a^{n-1})=\varphi(h, a^0, \dots, a^i, 1,
a^{i+1}, \dots, a^{n-1}); 0\le j\le n-1, \\
&t_n\varphi(h, a^0, \dots, a^n)=\varphi(h\ps{2},
S^{-1}(h\ps{1})(a^n), a^0, \dots, a^{n-1}).
\end{align}

We recall a cyclic map between  co-cyclic modules  is a linear map
which commutes with co-simplicial and cyclic operators.
\begin{proposition}
The following map is cyclic.
\begin{equation}
\Psi:C^n_\Hc(A)\longrightarrow C^n(A\rt \Hc),
\end{equation}
\begin{align*}
&\Psi(\varphi)( a^0\ot h^1, \dots, a^n\ot
h^n)=\\
&\varphi(h^0\ps{n+1}\dots h^{n-1}\ps{2}h^n, a^0, h^0\ps{1}(a^1),
\dots,  h^0\ps{n}\dots h^{n-1}\ps{1}(a^n)).
\end{align*}
\end{proposition}
\begin{proof}
We need to show that $\Psi$ commutes with cyclic structure of both
hand sides. The only non trivial ones are the commutation with
cyclic operators and the last face maps. We only show the former and
leave the rest to  the reader. Let $\tau_n$ be the cyclic operator
of the co-cyclic module $C^*(A\rt \Hc)$. Indeed,
\begin{align*}
&\tau_n(\Psi(\varphi))( a^0\ot h^1,\dots, a^n\ot
h^n)=\\
&\Psi(\varphi)(a^n\ot h^n,a^0\ot h^0 \dots, a^{n-1}\ot h^{n-1})=\\
&\varphi( h^n\ps{n+1}h^0\ps{n}\dots h^{n-2}\ps{2}h^{n-1}, a^n,
h^n\ps{1}(a^0), \dots, h^n\ps{n}h^0\ps{n-1}\dots
h^{n-2}\ps{1}(a^{n-1}))=\\
&\varphi(h^n\ps{n+1}\cdot(h^0\ps{n}\dots
h^{n-2}\ps{2}h^{n-1}h^n\ps{n+2}),a^n, h^n\ps{1}(a^0), \dots\\
&\dots, h^n\ps{n}h^0\ps{n-1}\dots h^{n-2}\ps{1}(a^{n-1}))=\\
&\varphi(h^0\ps{n}\dots
h^{n-2}\ps{2}h^{n-1}h^n\ps{n+2},S^{-1}(h^n\ps{n+1})\cdot (a^n,
h^n\ps{1}(a^0), \dots\\
 &\dots, h^n\ps{n}h^0\ps{n-1}\dots
h^{n-2}\ps{1}(a^{n-1})))= \\
&\varphi(h^0\ps{n}\dots
h^{n-2}\ps{2}h^{n-1}h^n\ps{2},S^{-1}(h^n\ps{1}) (a^n), a^0,
h^0\ps{1}(a^1), \dots\\
 &\dots, h^0\ps{n-1}\dots h^{n-2}\ps{1}(a^{n-1}))=\\
&\Psi(t\varphi)(a^0\ot h^0, \dots, a^n\ot h^n ).
\end{align*}
\end{proof}
We are ready now to determine the image of $\Psi$. The proof of the
following proposition is obvious.
\begin{proposition}
The map $\Psi$ lands in $C^*(A\rt \Hc; \Hc)$ provided we restrict
$\Psi$ to the normalized subcomplex of $C^*_\Hc(A)$.
\end{proposition}
\begin{theorem}\label{th}
The restriction of $\Psi$ to the normalized co-chains  is an
isomorphism and hence
\begin{equation}
HC^*_\Hc(A)\cong HC^*(A\rt\Hc; \Hc).
\end{equation}
\end{theorem}
\begin{proof}
Let us by abuse of notation denote the restriction of $\Psi$ on the
normalized complex again by $\Psi$. Using the above propositions, to
finish the proof we need the inverse of $\Psi$. Let $\Phi: C^*(A\rt
\Hc; \Hc)\ra \bar{C}^\ast_\Hc(A)$, where $\bar{C}$  denotes the
normalized complex,  be defined by:
$$(\Phi f)(h,a^0, \dots, a^n)= f(a^0\ot 1, \dots, a^{n-1}\ot 1, a^n\ot h).$$
First we need to show  that $\Phi$ is well-defined. To this end we
check that  $\Phi( f)$ is equivariant if $f$ is $\Hc$-constant. Let
$f$ be $\Hc$-constant. So,   $$bf(a^0\ot h^0,1\ot h,a^2\ot h^2,
\dots, a^n\ot h^n)=0.$$ We conclude that,
\begin{align}\label{yek}
&f(a^0\ot h^0h, a^2\ot h^2, \dots, a^n\ot h^n)=\\\notag
 &f(a^0\ot
h^0, h\ps{1}(a^1)\ot h\ps{2}h^1, a^2\ot h^2, \dots, a^n\ot h^n ).
\end{align}
 Similarly, for each $i\le n-1$,  we have the following identity,

\begin{align}\label{do}
&f(a^0\ot h^0, \dots,a^ih\ot h^i, \dots, a^n\ot h^n)=\\\notag
&f(a^0\ot h^0,\dots, h\ps{1}(a^{i+1})\ot  h\ps{2}h^{i+1}, a^2\ot
h^2, \dots, a^n\ot h^n )
\end{align}
And eventually we have,

\begin{align}\label{seh}
& f(a^0\ot h^0, \dots,a^{n-1}\ot h^{n-1}, a^n\ot h^nh)=\\\notag
&f(h\ps{1}(a^0)\ot h\ps{2}h^0, a^1\ot h^1, \dots, a^n\ot h^n ).
\end{align}
Now, by using \eqref{yek},\eqref{do}, \eqref{seh}, and Lemma
\ref{condition},  we show $\Phi f$ is equivariant,
\begin{align*}
& (\Phi f)(h,g\ps{1}(a^0),\dots, g\ps{n+1}(a^n))= \\
&f(g\ps{1}(a^0)\ot 1, \dots, g\ps{n}(a^{n-1})\ot 1,
g\ps{n+1}(a^n)\ot h)=\\
&f(g\ps{1}(a^0)\ot 1, \dots, g\ps{n}(a^{n-1})\ot 1,
g\ps{n+1}(a^n)\ot g\ps{n+2}S(g\ps{n+3})h)=\\
&f(g\ps{1}(a^0)\ot 1, \dots, g\ps{n}(a^{n-1})\ot g\ps{n+1},
a^n\ot S(g\ps{n+3})h)=\\
\vdots\\
&f(g\ps{1}(a^0)\ot g\ps{2}, a^1\ot 1, \dots, a^{n-1}\ot 1,
a^n\ot S(g\ps{3})h)=\\
&f(a^0\ot 1, \dots, a^{n-1}\ot 1, a^n\ot S(g\ps{2})hg\ps{1})=\Phi
f(S(g\ps{2})hg\ps{1}, a^0, \dots, a^n).
\end{align*}

The only thing left to finish the proof  is to show that $\Psi$ and
$\Phi$ are inverse to one another.  It is obvious that
$\Phi\Psi=id$. Again by using \eqref{yek}, \eqref{do}, and
\eqref{seh},  we show the $\Psi\Phi=id$.
\begin{align*}
&\Psi\Phi f(a^0\ot h^0, \dots, a^n\ot h^n)=\\
&\Phi f(h^0\ps{n+1}\dots h^{n-1}\ps{2}h^n, a^0, h^0\ps{1}(a^1),
\dots,  h^0\ps{n}\dots h^{n-1}\ps{1}(a^n))=\\
&f(a^0\ot 1, h^0\ps{1}(a^1)\ot 1, \dots,h^0\ps{n-1}\dots
h^{n-2}\ps{1}(a^{n-1})\ot 1,\\
&   h^0\ps{n}\dots h^{n-1}\ps{1}(a^n)\ot h^0\ps{n+1}\dots
h^{n-1}\ps{2}h^n)=\\
&f(a^0\ot h^0, a^1\ot 1, h^1\ps{1}(a^2)\ot 1 \dots, h^1\ps{n-2}\dots
h^{n-2}\ps{1}(a^{n-1})\ot 1,\\
&   h^1\ps{n-1}\dots h^{n-1}\ps{1}(a^n)\ot h^1\ps{n}\dots
h^{n-1}\ps{2}h^n)=\\
&f(a^0\ot h^0, a^1\ot h^1, a^2\ot 1, h^2\ps{1}(a^3)\ot 1,\dots,
h^2\ps{n-3}\dots
h^{n-2}\ps{1}(a^{n-1})\ot 1,\\
&   h^2\ps{n-2}\dots h^{n-1}\ps{1}(a^n)\ot h^2\ps{n-1}\dots
h^{n-1}\ps{2}h^n)=\\
\vdots\\
&f(a^0\ot h^0, \dots, a^n\ot h^n)
\end{align*}
\end{proof}
 It is shown  that  the equivariant cyclic cohomology and cyclic
cohomology of crossed product algebra coincide if $\Hc$ is
semisimple \cite{ak}.

\begin{corollary}
If $\Hc$ is  semisimple,  then $HC^*(A\rt \Hc; \Hc)\cong HC^*(A\rt
\Hc)$.
\end{corollary}

Now let $B$ be  a unital sub $\Hc$-module algebra of $A$. That is,
$B$ is a unital  subalgebra of $A$ and $h(b)\in B$, for any $h\in
\Hc$ and $b\in B$. Similarly with the same proof as of Theorem
\ref{th} one can show that
$$HC^*_\Hc(A;B)\cong HC^*(A\rt\Hc; B\rt H),$$
where the left hand side is the cyclic cohomology of subcomplex of
$B$-constant equivariant co-chains; i.e, $\varphi\in {\rm
Hom}_\Hc(\Hc\ot A^{\ot (n+1)}, \bbc)$, such that
$\varphi(h,a^0,\dots,a^n)$, and $b\varphi(h, a^0, \dots, a^{n+1})$
vanish if at least one of $a^j$ is in $ B$, for  $ j\ge 1.$

\vspace{2cm} E-mail address: bahram@math.ohio-state.edu
\end{document}